\documentclass[british]{article}
\usepackage{ae,aecompl}
\usepackage[T1]{fontenc}
\usepackage[latin9]{luainputenc}
\usepackage{amsmath}
\usepackage{amsthm}
\usepackage{amssymb}
\usepackage{showkeys}
\makeatletter
\newcommand{\lyxaddress}[1]{
	\par {\raggedright #1
	\vspace{1.4em}
	\noindent\par}
}
\theoremstyle{plain}
\newtheorem{thm}{\protect\theoremname}
\theoremstyle{remark}
\newtheorem{rem}[thm]{\protect\remarkname}
\theoremstyle{definition}

\theoremstyle{plain}
\newtheorem{prop}[thm]{\protect\propositionname}
\theoremstyle{plain}
\newtheorem{lem}[thm]{\protect\lemmaname}
\newtheorem{cor}[thm]{\protect\corollaryname}

\@ifundefined{date}{}{\date{}}

\usepackage{babel}
\addto\captionsbritish{\renewcommand{\definitionname}{Definition}}
\addto\captionsbritish{\renewcommand{\lemmaname}{Lemma}}
\addto\captionsbritish{\renewcommand{\propositionname}{Proposition}}
\addto\captionsbritish{\renewcommand{\remarkname}{Remark}}
\addto\captionsbritish{\renewcommand{\theoremname}{Theorem}}
\addto\captionsenglish{\renewcommand{\definitionname}{Definition}}
\addto\captionsenglish{\renewcommand{\lemmaname}{Lemma}}
\addto\captionsenglish{\renewcommand{\propositionname}{Proposition}}
\addto\captionsenglish{\renewcommand{\remarkname}{Remark}}
\addto\captionsenglish{\renewcommand{\theoremname}{Theorem}}
\providecommand{\definitionname}{Definition}
\providecommand{\lemmaname}{Lemma}
\providecommand{\propositionname}{Proposition}
\providecommand{\remarkname}{Remark}
\providecommand{\theoremname}{Theorem}
\providecommand{\corollaryname}{Corollary}

\makeatother

\usepackage{babel}
\addto\captionsbritish{\renewcommand{\definitionname}{Definition}}
\addto\captionsbritish{\renewcommand{\lemmaname}{Lemma}}
\addto\captionsbritish{\renewcommand{\propositionname}{Proposition}}
\addto\captionsbritish{\renewcommand{\remarkname}{Remark}}
\addto\captionsbritish{\renewcommand{\theoremname}{Theorem}}
\addto\captionsenglish{\renewcommand{\definitionname}{Definition}}
\addto\captionsenglish{\renewcommand{\lemmaname}{Lemma}}
\addto\captionsenglish{\renewcommand{\propositionname}{Proposition}}
\addto\captionsenglish{\renewcommand{\remarkname}{Remark}}
\addto\captionsenglish{\renewcommand{\theoremname}{Theorem}}
\providecommand{\definitionname}{Definition}
\providecommand{\lemmaname}{Lemma}
\providecommand{\propositionname}{Proposition}
\providecommand{\remarkname}{Remark}
\providecommand{\theoremname}{Theorem}

\newcommand{\IN}{\mathbb{N}}
\newcommand{\IR}{\mathbb{R}}
\newcommand{\IE}{\mathbb{E}}
\newcommand{\IP}{\mathbb{P}}
\newcommand{\ui}{\underline{i}}
\newcommand{\uj}{\underline{j}}
\newcommand{\Ne}{\mathcal{N}_{1}}
\newcommand{\Nz}{\mathcal{N}_{2}}

\begin{document}
\title{Two Groups in a Curie-Weiss Model with Heterogeneous Coupling}
\author{Werner Kirsch and Gabor Toth}
\maketitle

\lyxaddress{Fakult\"at f\"ur Mathematik und Informatik \\
FernUniversit\"at Hagen\\
D-58084 Hagen\\
Germany\\
email: werner.kirsch@fernuni-hagen.de\qquad{}gabor.toth@fernuni-hagen.de }
\begin{abstract}
\noindent We discuss a Curie-Weiss model with two groups with different
coupling constants within and between groups. For the total magnetisations
in each group, we show bivariate laws of large numbers and a central
limit theorem which is valid in the high temperature regime. In the
critical regime, the total magnetisation normalised by $N^{3/4}$
converges to a non-trivial distribution which is not Gaussian, just
as in the single-group Curie-Weiss model. Finally, we prove a kind of a
`law of large numbers' in the low temperature
regime, more precisely we prove that the empirical magnetisation converges
in distribution to a mixture of two Dirac measures.

\smallskip{}


\smallskip{}

\end{abstract}

\section{Introduction}

The Curie-Weiss model is probably the easiest model of magnetism which
shows a phase transition between a paramagnetic and a ferromagnetic
phase. In this model the spins can take values in $\{-1,1\}$ (or
up/down), each spin interacts with all the others in the same way.
More precisely, for finitely many spins $X:=(X_{1},X_{2},\ldots,X_{N})\in\{-1,1\}$
the energy of the spins is given by
\begin{align}
H~=~H(X_{1},\ldots,X_{N})~:=~-\frac{1}{2N}\,\big(\sum_{j=1}^{N}\,X_{j}\big)^{2}\,.
\end{align}

The `Gibbs measure' or `canonical ensemble' with coupling constant ($\approx$ inverse temperature) $J_{0}\geq0$ the
probability of a spin configuration is given by
\begin{align}
\mathbb{P}\big(X_{1}=x_{1},\ldots,X_{N}=x_{N}\big)~:=~Z^{-1}\;e^{-J_{0}H(x_{1},\ldots,x_{N})}\label{eq:PJ_0}
\end{align}
where $x_{i}\in\{-1,1\}$ and $Z$ is a normalisation constant which
depends on $N$ and $J_{0}$.

The quantity
\begin{align}
S_{N}~=~\sum_{j=1}^{N}X_{j}
\end{align}
is called the (total) magnetisation. It is well known (see e.\,g. \marginpar{more citations}
Ellis \cite{Ellis} or \cite{MM}) that the Curie-Weiss model has
a phase transition at $J_{0}=1$ in the following sense
\begin{align}
\frac{1}{N}\,S_{N}~\Longrightarrow~\frac{1}{2}\,(\delta_{-m(J_{0})}+\delta_{m(J_{0})})\label{eq:lln}
\end{align}
where $\Rightarrow$ denotes convergence in distribution, $\delta_{x}$
the Dirac measure in $x$.

For $J_{0}\leq1$ we have $m(J_{0})=0$ which is the unique solution
of
\begin{align}
\tanh(J_{0}x)=x\label{eq:mJ_0}
\end{align}
for this case.

If $J_{0}>1$ equation (\ref{eq:mJ_0}) has exactly three solutions
and $m(J_{0})$ is the unique positive one.

Equation (\ref{eq:lln}) is a substitute for the law of large numbers
for i.i.d. random variables.

Moreover, for $J_{0}<1$ there is a central limit theorem, i.\,e.
\begin{align}
\frac{1}{\sqrt{N}}\,S_{N}~\Longrightarrow~\mathcal{N}(0,\frac{1}{1-J_{0}})
\end{align}

There is a huge amount of literature on the Curie-Weiss model. We can just mention a few papers here. 
The Curie-Weiss model is also called the Husimi-Temperley model. It
was first introduced by Husimi \cite{Husimi} and Temperley \cite{Temperley}.
Subsequently, it was discussed by Kac \cite{Kac} and Ellis and Newman \cite{EN1,EN2}. It was also treated in the textbooks Thompson \cite{Thompson},
and Ellis \cite{Ellis}.

More recently, the Curie-Weiss model has
been used in the context of social and political interactions. See
e.g. \cite{CGh, GBC, Penrose}.

In this paper we consider two groups of Curie-Weiss spins
$X~=~(X_{1},\ldots,X_{{N}_{1}})$ and $Y~=~(Y_{1},\ldots,Y_{{N}_{2}})$
with $N:=N_{1}+{N}_{2}$. The spins $X$ and $Y$ are Curie-Weiss spins
with coupling constant $J_{1}$ and $J_{2}$ respectively, in addition there is a Curie-Weiss-type interaction between the $X_{i}$
and the $Y_{j}$ with coupling constant $\bar{J}$.

We set
\begin{align} J:= \begin{pmatrix}J_{1}&\bar{J}\\\bar{J}&J_{2}\end{pmatrix} \end{align}
and assume that $J_{1},J_{2}, \bar{J}>0$ and
\begin{align}\label{eq:delta} \Delta:=J_{1}J_{2}-\bar{J}^{2}>0\,, \end{align}
so that the matrix $J$ is positive definite. Loosely speaking, this conditions ensures that the interaction within groups
dominates the interaction between the groups.

The energy function is given by
\begin{align}
H=H_{J}(X,Y):=-\frac{1}{2N}\left[J_{1}\big(\sum_{j=1}^{{N}_{1}}X_{j}\big)^{2}+J_{2}\big(\sum_{j=1}^{{N}_{2}}Y_{j}\big)^{2}
+2\bar{J}\sum_{i=1}^{{N}_{1}}\sum_{j=1}^{{N}_{2}}X_{i}Y_{j}\right].
\end{align}

We denote the Gibbs measure associated with $H_{J}(X,Y)$ by $\IP_{J}$ (sometimes abbreviated by $\IP$) defined by
\begin{align}
   \IP_{J}\big(A\big)~:=~Z^{-1}\;\sum_{(X,Y)\in A} e^{-H_{J}(X,Y)}
\end{align}
where $Z$ is a normalizing constant which makes $\IP_{J}$ a probability measure. The corresponding expectation is called $\IE_{J}$, sometimes abbreviated $\IE$.

By sending $N$ to infinity we mean that both ${N}_{1}$ and
${N}_{2}$ tend to infinity. We set
\begin{equation}
\alpha_{1}:=\lim\frac{N_{1}}{N},\qquad\qquad\alpha_{2}:=\lim_{N\to\infty}\frac{N_{2}}{N}~=~1-\alpha_{1}
\end{equation}
and assume that these limits exist and $0<\alpha_{1}<1$.

In this paper we consider the asymptotic behaviour of the two-dimensional random
variables
\begin{equation}\label{mean}
\big(\frac{1}{{N_{1}}^{\gamma}}\sum_{i=1}^{N_{1}}\,X_{i}\,,\,\frac{1}{{N_{2}}^{\gamma}}\sum_{j=1}^{N_{2}}\,Y_{j}\,\big)
\end{equation}
where $\gamma=1, \frac{1}{2}$ or $\frac{3}{4}$ depending on the parameters of the model, namely $J_{1}, J_{2}, \bar{J}$ and $\alpha_{1}, \alpha_{2}$.

Our assumptions on $J_{1}, J_{2}$ and $\bar{J}$ exclude a few `borderline' cases. If $\bar{J}=0$ the two groups are independent
of each other and can therefore be treated as independent single group Curie-Weiss spins. This is also the case if $J_{1}=0$
or $J_{2}=0$ as this implies $\bar{J}=0$ by assumption \eqref{eq:delta}.

Condition \ref{eq:delta} also excludes the (interesting) case $J=\begin{pmatrix}J_{0}&J_{0}\\ J_{0}& J_{0}\end{pmatrix}$. This case which we call the homogeneous one
requires a somewhat different technique. It is treated in \cite{KT1}.

Another borderline case is given by $\alpha_{1}=0 $ or $\alpha_{2}=0 $ (assuming still that both $N_{1}$ and $N_{2}$ tend to infinity). We may even consider
the following extension of our model: The groups may consist of $\tilde{N}_{i}\approx \rho_{i}N$
but the averages in \eqref{mean} are taken over $N_{i}\approx \alpha_{i}N\leq \tilde{N}_{i}$.

These cases can be treated by the techniques of this paper as well. With the obvious changes, the results and their proofs remain
valid for these extensions.
In order to avoid
a notational overkill we stick to the stronger assumptions made.

\subsection{High Temperature Regime}

The parameter space of this model is

\[
\Phi:=\{(J_{1},J_{2},\bar{J},\alpha_{1},\alpha_{2})\in(0,\infty)^{5}|J_{1}J_{2}>\bar{J}^{2},\alpha_{1}+\alpha_{2}=1\}.
\]

For the single-group model, the high temperature regime is quite simply
expressed by the single condition $J_{0}<1$. For two groups with
a heterogeneous coupling matrix, we have a very different situation:
each within-group coupling constant $J_{\nu}$ has to be small in
relation to the reciprocal of the group's size. Once the within-group
couplings have been chosen, the between-groups coupling has to be
small, too. How small depends on how close the other two couplings
are to the reciprocals of the group sizes. If the within-group couplings
are very small, that leaves more leeway for the between-groups coupling
to be larger.

We shall assume that the interactions satisfy

\begin{align}
J_{1} & \;<\;\frac{1}{\alpha_{1}},\label{eq:condJ1-1}\\
J_{2} & \;<\;\frac{1}{\alpha_{2}},\label{eq:condJ2-1}\\
\bar{J}^{2} & \;<\;(\frac{1}{\alpha_{1}}-J_{1})(\frac{1}{\alpha_{2}}-J_{2}).\label{eq:condJ3-1}
\end{align}
and refer to these conditions as the `high temperature regime', and
we shall also refer to the subset $\Phi_{h}$ of $\Phi$ where these
conditions hold by the same name. Note that if we use the symbol $\alpha$
for the diagonal $2\times2$ matrix with entries $\alpha_{1}$ and
$\alpha_{2}$, we can formulate these conditions equivalently in matrix
form: the matrix

\[
J^{-1}-\alpha
\]
is positive definite if and only if we are in the high temperature
regime (see Proposition \ref{prop:equivalence}).

We prove a `law of large numbers'.
\begin{thm}
\label{LLN} In the high temperature regime, we have

\[
\big(\frac{1}{N_{1}}\sum_{i=1}^{N_{1}}X_{i},\frac{1}{N_{2}}\sum_{j=1}^{N_{2}}Y_{j}\big)~\underset{N\to\infty}{\Longrightarrow}(0,0).
\]
\end{thm}

Above '$\Longrightarrow$' denotes convergence in distribution of
the $2$-dimensional random variable on the left hand side.

We also have a `central limit theorem'. Using $\alpha=\begin{pmatrix}
                                                                              \alpha_{1}& 0\\ 0 & \alpha_{2}
                                                                           \end{pmatrix}$
we define the matrix
\begin{align}\label{eq:covariance}
   C~=~1+\Big(\alpha^{-1/2}\,J^{-1}\,\alpha^{-1/2}\;-\;1\Big)^{-1}
\end{align}
  where $1$ denotes the identity matrix.

\begin{thm}
\label{CLT} In the high temperature regime, we have
\begin{equation}
(\frac{1}{\sqrt{N_{1}}}\sum_{i=1}^{N_{1}}X_{i},\frac{1}{\sqrt{N_{2}}}\sum_{j=1}^{N_{2}}Y_{j})~\underset{N\to\infty}{\Longrightarrow}~\mathcal{N}\big((0,0),C),\label{eq:CLT}
\end{equation}
The covariance matrix $C$ (as in \eqref{eq:covariance}) is given by
\begin{align}
C & =\frac{1}{(1-\alpha_{1}J_{1})(1-\alpha_{2}J_{2})-\alpha_{1}\alpha_{2}\bar{J}^{2}}\left[\begin{array}{cc}
1-\alpha_{2}J_{2} & \sqrt{\alpha_{1}\alpha_{2}}\bar{J}\\
\sqrt{\alpha_{1}\alpha_{2}}\bar{J} & 1-\alpha_{1}J_{1}
\end{array}\right]\label{eq:covariance1}
\end{align}
\end{thm}
\begin{rem}
   Theorem \ref{CLT} implies that also expressions like $\frac{1}{\sqrt{N_{1}}}\sum_{i=1}^{N_{1}}X_{i}\pm\frac{1}{\sqrt{N_{2}}}\sum_{j=1}^{N_{2}}Y_{j}$
are asymptotically Gaussian distributed.
\end{rem}

\subsection{Critical Regime}

The critical regime is where $\big(\sum_{i=1}^{N_{1}}X_{i},\sum_{j=1}^{N_{2}}Y_{j}\big)$
abruptly changes behaviour. In the single-group model, this occurs
at $J_{0}=1$. For two groups with a heterogeneous coupling matrix,
in the critical regime, each within-group coupling constant $J_{\nu}$
has to be small in relation to the reciprocal of the group's size.
Once the within-group couplings have been chosen, the between-groups
coupling has to have an exact magnitude, which is larger than in the
high temperature regime:
\begin{align}
J_{1}~ & <~\frac{1}{\alpha_{1}},\label{eq:condJ1-1-1}\\
J_{2}~ & <~\frac{1}{\alpha_{2}},\label{eq:condJ2-1-1}\\
\bar{J}^{2}~ & =~(\frac{1}{\alpha_{1}}-J_{1})(\frac{1}{\alpha_{2}}-J_{2}).\label{eq:condJ3-1-1}
\end{align}

We shall call the subset of $\Phi$ where these conditions hold $\Phi_{c}$
and we can also formulate these conditions equivalently in matrix
form: the matrix

\[
J^{-1}-\alpha
\]
is singular and has positive diagonal entries if and only if we are
in the critical regime.

We also note that if $J_{1}=\frac{1}{\alpha_{1}}$ or
$J_{2}=\frac{1}{\alpha_{2}}$ then \eqref{eq:condJ3-1-1} implies $\bar{J}=0$ hence the
two groups are independent of each other and can be treated as in the single-group case.

In the critical regime we consider here \big(i.\,e. for (\ref{eq:condJ1-1-1})--(\ref{eq:condJ3-1-1})\big)
the law of large numbers, Theorem \ref{LLN}, still holds,  but the central limit theorem, Theorem \ref{CLT},
has to be replaced by a theorem describing the asymptotic behaviour
of
\begin{align}
T_{N}~=~\Big(\frac{1}{{N_{1}}^{3/4}}\,S_{N_{1}}^{(1)},\frac{1}{{N_{2}}^{3/4}}\,S_{N_{2}}^{(2)}\Big)\,.\label{eq:cltcrit}
\end{align}
This sequence $T_{N}$ converges in distribution but not to a normal
distribution. We state the moments of the limiting measure in Theorem
\ref{Fluctuations}.

The critical regime results are:
\begin{thm}
\label{LLN-1} In the critical regime, we have

\[
\big(\frac{1}{N_{1}}\sum_{i=1}^{N_{1}}X_{i},\frac{1}{N_{2}}\sum_{j=1}^{N_{2}}Y_{j}\big)~\underset{N\to\infty}{\Longrightarrow}\delta_{(0,0)}.
\]
\end{thm}

If we choose as normalising factors $N_{\nu}^{\frac{3}{4}}$ instead
of $N_{\nu}$, then we obtain
\begin{thm}
\label{Fluctuations} In the critical regime, the random variables

\[
(\frac{1}{N_{1}^{3/4}}\sum_{i=1}^{N_{1}}X_{i},\frac{1}{N_{2}^{3/4}}\sum_{j=1}^{N_{2}}Y_{j})
\]

converge in distribution to a measure $\mu$ (on $\mathbb{R}^{2}$)
with moments

\begin{align}
 & m_{K,Q}~:=~\int x^{K}y^{Q}\;d\mu(x,y)\notag\\
~ & =~\left[\frac{12}{\alpha_{1}(L_{2}-\alpha_{2})^{2}+\alpha_{2}(L_{1}-\alpha_{1})^{2}}\right]^{\frac{K+Q}{4}}(L_{1}-\alpha_{1})^{\frac{Q}{2}}(L_{2}-\alpha_{2})^{\frac{K}{2}}\cdot\notag\\
 & \quad\cdot\frac{\Gamma(\frac{K+Q+1}{4})}{\Gamma(\frac{1}{4})}\alpha_{1}^{\frac{K}{4}}\alpha_{2}^{\frac{Q}{4}}.
\end{align}
\end{thm}
where the matrix $L$ is defined by $L= \begin{pmatrix}
                                          L_{1}&-\bar{L}\\-\bar{L}&L_{2}
                                       \end{pmatrix}=J^{-1}$

\subsection{Low Temperature Regime}

In the single-group model, the low temperature regime is characterised
by the inequality $J_{0}>1$. The magnetization $\frac{1}{N}\sum X_{i}$ converges
in distribution to the measure $\frac{1}{2}(\delta_{m^{*}}+\delta_{-m^{*}})$ where $m^{*}$
is the unique positive solution of the equation $m=\tanh(J_{0}m)$. We regard this fact as a (substitute
for the) law of large numbers.

In the case of two groups we define the low temperature
regime to be the complement $\Phi_{l}:=\Phi\backslash(\Phi_{h}\cup\Phi_{c})$
in the parameter space. We have a similar `law of large numbers' in this case.
\begin{thm}
\label{thm:SpecialCase}In the low temperature regime,
there are exactly two nonzero solutions $m^{*}=(m^{*}_{1},m^{*}_{2}) $ and $-m^{*}$ of the system
\begin{align}
   m_{1}~&=~\tanh(J_{1}\alpha_{1}m_{1}-\bar{J}\alpha_{2}m_{2})\\
   \text{and}\qquad m_{2}~&=~\tanh(J_{2}\alpha_{2}m_{2}-\bar{J}\alpha_{1}m_{1})\,.
\end{align}
We have
\begin{align}
   \big(\frac{1}{N_{1}}\sum_{i=1}^{N_{1}}X_{i},\frac{1}{N_{2}}\sum_{j=1}^{N_{2}}Y_{j}\big)~\underset{N\to\infty}{\Longrightarrow}\frac{1}{2}\big(\delta_{-m^{*}}
   +\delta_{m^{*}}\big)\,.
\end{align}

Moreover, we may assume $m^{*}_{1}>0$ and $m^{*}_{2}\,>0$.
\end{thm}

Obviously, in the low temperature case, there can be no central limit theorem in the sense that $(\frac{1}{N_{1}^{\gamma}}\sum X_{i},\frac{1}{N_{2}^{\gamma}}\sum Y_{i})$ converges to a nontrivial limit measure. However, we conjecture that there is
a `conditional version' of a central limit theorem.

\textbf{Acknowledgment:} While finishing this paper we became aware of the papers \cite{FM}
and \cite{FC} which contain the above results as special cases. The
methods used by those authors is very different from ours. We are
grateful to Francesca Collet for drawing our attention to the papers
\cite{FM} and \cite{FC}.

We would also like to thank Matthias L\"owe and Kristina Schubert \cite{LS} as well as an unnamed referee for
valuable comments which in our opinion improved this paper considerably.

\section{A rough sketch of the proofs}

Our proofs are based on the method of moments, the basis of which is the following well known Theorem (see e.~g. \cite{Breiman}).
\begin{thm}\label{thm:moments}
   Suppose that $\mu_{N}$ and $\mu $ are probability measure on $\IR^{d}$ for which all moments are finite and assume that $\mu $ is determined by its moments $m_{\underline{k}}(\mu) $.
   If $m_{\underline{k}}(\mu_{N}) \to m_{\underline{k}}(\mu)$ for all $\underline{k}\in\IN^{d}$ as $N\to\infty$ then the measures $\mu_{N} $ converge weakly to $\mu$.
\end{thm}
It is also known that all (multidimensional) normal distributions are determined by their moments (see e.~g. \cite{KS}).

Thus we will consider suitably normalized moments of the form
\begin{align}\label{eq:moments}
   M_{K,Q}~:=&~\IE_{J}\Big(\big(\sum_{i=1}^{N_{1}}X_{i}\big)^{K}\big(\sum_{j=1}^{N_{2}}Y_{j}\big)^{Q}\Big)\notag\\
   =& \sum_{i_{1},i_{2},\ldots,i_{K}}\;\sum_{j_{1},j_{2},\ldots,j_{Q}}\;\IE_{J}\Big(X_{i_{1}}\cdot\ldots\cdot X_{i_{K}}\,
   \cdot\,Y_{j_{1}}\cdot\ldots\cdot Y_{j_{Q}}\Big)\\
   =&~\sum_{\ui\in\Ne^{K}}\;\sum_{\uj\in\Nz^{Q}}\;\IE_{J}\Big(X(\ui)\,
   \cdot\,Y(\uj)\Big)\notag
\end{align}
where $\ui=(i_{1},\ldots,i_{K})$, $\uj=(j_{1},\ldots,j_{Q})$, $\Ne=\{1,2,\ldots,N_{1}\}$, $\Nz=\{1,2,\ldots,N_{2}\}$ and $X(\ui)=\prod_{\nu=1}^{K}X_{i_{\nu}}$.

Since ${X_{i}}^{2}={Y_{j}}^{2}=1$ and due to exchangeability we have
\begin{align}
   &\IE_{J}\Big(X_{i_{1}}\cdot\ldots\cdot X_{i_{K}}\,
   \cdot\,Y_{j_{1}}\cdot\ldots\cdot Y_{j_{Q}}\Big)\notag\\
   =~&\IE_{J}\Big(X_{1}\cdot\ldots\cdot X_{\tilde{K}}\,
   \cdot\,Y_{1}\cdot\ldots\cdot Y_{\tilde{Q}}\Big)\label{eq:simpmoment}
\end{align}
where $\tilde{K}$ (resp. $\tilde{Q}$) is the number of $i_{k}$ (resp. $j_{\ell}$) which occur an odd number of times.

In order to evaluate the moments in \eqref{eq:moments}, we need to estimate correlations as in \eqref{eq:simpmoment}.

In section \ref{sec:correl} we will prove asymptotic estimates for the correlations. For example in section \ref{sub:corrhigh} show that \emph{in the high temperature regime}
\begin{align}\label{eq:correl}
   \Big|\IE_{J}\Big(X_{1}\cdot\ldots\cdot X_{{K}}\,
   \cdot\,Y_{1}\cdot\ldots\cdot Y_{{Q}}\Big)\Big|~
   \leq~c_{K,Q}\;\frac{1}{N^{(K+Q)/2}}\;
\end{align}
where $c_{K,Q}$ depends on the matrix $J$ and the numbers $\alpha_{1}$ and $\alpha_{2}$ but not on $N$. Note that in the case of independent random variables, i.~e. if
$J=\begin{pmatrix} 0 & 0\\ 0 & 0\\ \end{pmatrix}$ the correlation \eqref{eq:correl} is zero unless $K=Q=0$.

In a second step of the proof we need a `bookkeeping' method, to keep track of the variety of terms in the sum \eqref{eq:moments}.
We have to count the number of simple, double, triple, etc. occurrences of the $X_{i}$ and $Y_{j}$ in \eqref{eq:moments}.
We start the discussion of this topic with Lemma \ref{thm:comb1}. A more precise discussion is given in subsection \ref{subs:comb}.
In section \ref{sec:lln} we prove the laws of large numbers, Theorems \ref{LLN} and \ref{LLN-1}. The central limit theorem
is proved in \ref{sub:prclt} combining the results on the correlations and the bookkeeping method.

\section{Computing Expectations}\label{sec:correl}
In this section we compute expectations (=correlations) of the form
\begin{align}\label{eq:correlations}
   \IE_{J}\Big(X_{1}\cdot X_{2}\cdot\ldots\cdot X_{K}\,\cdot\,Y_{1}\cdot Y_{2}\cdot\ldots\cdot Y_{Q}\Big)
\end{align}
asymptotically for the three regimes of $J$.
\subsection{A two-dimensional Hubbard-Stratonovich transform}
For any configuration of the spins
\begin{align}
(X,Y)~=~(X_{1},X_{2},\ldots,X_{N_{1}},Y_{1},Y_{2},\ldots,Y_{N_{2}})
\end{align}
we set
\begin{align}
   S_{1}~=~\sum_{i=1}^{N_{1}}X_{i}&\qquad S_{2}=\sum_{j=1}^{N_{2}}Y_{j}\\
   \text{and}\qquad&S=
   \begin{pmatrix}
      S_{1}\\ S_{2}
   \end{pmatrix}
\end{align}
and define the function
\begin{align*}
h(S_{1},S_{2}) & =\frac{1}{2N}\,(S_{1},S_{2})J\left(\begin{array}{c}
S_{1}\\
S_{2}
\end{array}\right)\\
 & =\frac{1}{2N}\,S'JS,
\end{align*}

For a symmetric positive definite $2\times2$ matrix $A$ and a point
$x_{0}\in\mathbb{R}^{2}$ we can use the following equality to express
a value of the exponential function as an integral:

\[
e^{\frac{x_{0}'Ax_{0}}{2}}=\frac{\sqrt{\det A}}{2\pi}\int_{\mathbb{R}^{2}}e^{-\frac{x'Ax}{2}}e^{-x'Ax_{0}}\mathrm{d}x.
\]

According to this equality,

\[
e^{h(S_{1},S_{2})}=e^{\frac{S'JS}{2N}}=c\int_{\mathbb{R}^{2}}e^{-\frac{1}{2N}x'Jx}e^{-\frac{1}{N}x'JS}\mathrm{d}x,
\]
where $c=\frac{\sqrt{\det J}}{2\pi}$.

We define the inverse matrix

\begin{equation}
L=\left[\begin{array}{cc}
L_{1} & -\bar{L}\\
-\bar{L} & L_{2}
\end{array}\right]=\frac{1}{J_{1}J_{2}-\bar{J}^{2}}\left[\begin{array}{cc}
J_{2} & -\bar{J}\\
-\bar{J} & J_{1}
\end{array}\right]=J^{-1}.\label{eq:defL}
\end{equation}

Switching variables $y=\frac{1}{N}Jx$ we obtain

\begin{align}\label{eq:HS}
    e^{h(S_{1},S_{2})}=c' \int_{\mathbb{R}^{2}}e^{-\frac{N}{2}y\cdot Ly}e^{S\cdot y}\;\mathrm{d}^{2}y,
\end{align}
where $c'$ is a term that depends on the matrix $L$ and on $N$.
Equation \eqref{eq:HS} is our two-dimensional version of the Hubbard-Stratonovich transform.

Summing over all $(X,Y)\in \{-1,+1\}^{N}$ we obtain
\begin{align}
   \sum_{X,Y}\;e^{S\cdot y}~&=~(e^{+y_{1}}+e^{-y_{1}})^{N_{1}}\;\cdot\;(e^{+y_{2}}+e^{-y_{2}})^{N_{2}}\notag\\
   ~&=~2^{N}\,\cosh^{N_{1}}(y_{1})\,\cdot\,\cosh^{N_{2}}(y_{2})
\end{align}
and
\begin{align}
   &\sum_{X,Y}\; X_{1}\cdot\ldots X_{K}\,\cdot\,Y_{1}\cdot\ldots\cdot Y_{Q}\;e^{S\cdot y}\notag\\
   ~=~&\frac{(e^{+y_{1}}-e^{-y_{1}})^{K}}{(e^{+y_{1}}+e^{-y_{1}})^{K}}\,(e^{+y_{1}}+e^{-y_{1}})^{N_{1}}\;\cdot\;\frac{(e^{+y_{2}}-e^{-y_{2}})^{Q}}{(e^{+y_{2}}+e^{-y_{2}})^{Q}}\,(e^{+y_{2}}+e^{-y_{2}})^{N_{2}}\notag\\
   ~=~&2^{N}\,\tanh^{K}(y_{1})\,\tanh^{Q}(y_{2})\;\cosh^{N_{1}}(y_{1})\,\cosh^{N_{2}}(y_{2})
\end{align}

Consequently we have
\begin{align}
   &\sum_{X,Y}\; X_{1}\cdot\ldots X_{K}\,\cdot\,Y_{1}\cdot\ldots\cdot Y_{Q}\;e^{-H_{J}(X,Y)}\notag\\
     =~&c\,\int e^{-N(1/2y\cdot Ly-N_{1}/N\ln\cosh y_{1}-N_{2}/N\ln\cosh y_{2})}\tanh^{K}y_{1}\tanh^{Q}y_{2}\,\mathrm{d}^{2}y\notag\\
     =~&c\,\int e^{-N\,F_{J}(y)}\;\tanh^{K}y_{1}\,\tanh^{Q}y_{2}\,\mathrm{d}^{2}y.
\end{align}
where
\begin{align}
   F_{J}(y)~:=~\frac{1}{2}L_{1}y_{1}^{2}+\frac{1}{2}L_{2}y_{2}^{2}-\bar{L}y_{1}y_{2}-\alpha_{1}\ln\cosh y_{1}-\alpha_{2}\ln\cosh y_{2}.
\end{align}

Let us define
\begin{align}\label{eq:defZ}
   Z_{J}(K,Q)~:=~\int e^{-N F_{J}(y)}\;\tanh^{K}y_{1}\,\tanh^{Q}y_{2}\,\mathrm{d}^{2}y.
\end{align}
then
\begin{align}
   \IE_{J}\Big( X_{1}\cdot X_{2}\cdot\ldots\cdot X_{K}\,\cdot\,Y_{1}\cdot Y_{2}\cdot\ldots\cdot Y_{Q}\Big)~=~\frac{\,Z_{J}(K,Q)\,}{Z_{J}(0,0)}
\end{align}

Thus if we can compute $Z_{J}(K,Q)$ asymptotically we will be able to compute the correlations \eqref{eq:correlations}.

\subsection{Extrema of the function $F$}
\subsubsection{High Temperature Regime}
We are going to apply the Laplace method to evaluate the quantities $Z_{J}(K,Q)$. In order to do so we need to determine the minima
of the function

\begin{equation}
F(y_{1},y_{2})=\frac{1}{2}L_{1}y_{1}^{2}+\frac{1}{2}L_{2}y_{2}^{2}-\bar{L}y_{1}y_{2}-\alpha_{1}\ln\cosh y_{1}-\alpha_{2}\ln\cosh y_{2}.\label{eq:def_F}
\end{equation}

\begin{prop}
\label{prop:pos-def}If
\begin{align}
L_{1} & >\alpha_{1},\label{eq:posdef1}\\
(L_{1}-\alpha_{1})(L_{2}-\alpha_{2}) & >\bar{L}^{2},\label{eq:posdef2}
\end{align}

then the function $F$ has a unique minimum at $(0,0)$.

$F$ has strictly positive definite Hessian
\begin{align*}
  H=\left[\begin{array}{cc}
F_{11} & F_{12}\\
F_{21} & F_{22}
\end{array}\right]
\end{align*}

and is therefore strictly convex.
\end{prop}

We used $F_{ij}$ to denote the partial derivative of $F$ with respect to $y_{i}$ and $y_{j}$.

\begin{rem}
The conditions (\ref{eq:posdef1}) and (\ref{eq:posdef2}) are equivalent
to the high temperature regime, as we shall show in proposition \ref{prop:equivalence}.
\end{rem}

\begin{proof}
We take derivatives with respect to both variables

\begin{align}
F_{1}(y_{1},y_{2}) & =L_{1}y_{1}-\bar{L}y_{2}-\alpha_{1}\tanh y_{1},\label{eq:FOC1a}\\
F_{2}(y_{1},y_{2}) & =L_{2}y_{2}-\bar{L}y_{1}-\alpha_{2}\tanh y_{2},\label{eq:FOC2a}\\
F_{11}(y_{1},y_{2}) & =L_{1}-\frac{\alpha_{1}}{\cosh^{2}y_{1}},\nonumber \\
F_{22}(y_{1},y_{2}) & =L_{2}-\frac{\alpha_{2}}{\cosh^{2}y_{2}}.\nonumber
\end{align}

The Hessian matrix of $F$ is

\[
H=\left[\begin{array}{cc}
F_{11} & F_{12}\\
F_{21} & F_{22}
\end{array}\right]=\left[\begin{array}{cc}
L_{1}-\frac{\alpha_{1}}{\cosh^{2}y_{1}} & -\bar{L}\\
-\bar{L} & L_{2}-\frac{\alpha_{2}}{\cosh^{2}y_{2}}
\end{array}\right].
\]

One solution to the first order conditions (\ref{eq:FOC1a}) and (\ref{eq:FOC2a})
is $y_{1}=y_{2}=0$.

The matrix $H$ is positive definite at the origin if and only if
(\ref{eq:posdef1}) and (\ref{eq:posdef2}) hold. Hence there is a
local minimum at the origin. If the Hessian matrix at the origin is
positive definite, it is also positive definite at any other point
due to $\cosh|s|>\cosh|t|$ for all $|s|>|t|$. Thus $F$ is strictly
convex and it follows that the minimum is unique and global.
\end{proof}
\begin{lem}\label{lem:det}
\begin{align}
   {\rm sgn}\; \Big(\det\begin{pmatrix}
                     \frac{1}{\alpha_{1}}-J_{1} & \bar{J}\\ \bar{J}&  \frac{1}{\alpha_{2}}-J_{2}
                  \end{pmatrix}\Big)
                  ~=~{\rm sgn}\; \Big(\det\begin{pmatrix}
                     L_{1}-\alpha_{1} & -\bar{L}\\ -\bar{L}& L_{2}-\alpha_{2}
                  \end{pmatrix}\Big)
                  \end{align}
\end{lem}
We used the notation
\begin{align*}
   {\rm sgn}\,(x)~=~\left\{
                    \begin{array}{ll}
                      1, & \hbox{x>0;} \\
                      0, & \hbox{x=0;} \\
                      -1, & \hbox{x<0.}
                    \end{array}
                  \right.
\end{align*}
\begin{proof}
   We write
\begin{align*}
   \alpha^{-1}-J~=~J\,(J^{-1}-\alpha)\,\alpha^{-1}
\end{align*}
Since $\det(J),\det(\alpha)>0$ the assertion follows.
\end{proof}
\begin{prop}
\label{prop:equivalence}The conditions on the Hessian matrix H given
in Proposition \ref{prop:pos-def} are equivalent to the following
conditions on the coupling matrix $J$:

\begin{align}
J_{1} & <\frac{1}{\alpha_{1}},\label{eq:condJ1}\\
J_{2} & <\frac{1}{\alpha_{2}},\label{eq:condJ2}\\
\bar{J}^{2} & <(\frac{1}{\alpha_{1}}-J_{1})(\frac{1}{\alpha_{2}}-J_{2}).\label{eq:condJ3}
\end{align}
\end{prop}

\begin{proof}
The equivalence of (\ref{eq:posdef2}) and (\ref{eq:condJ3}) is the contents of Lemma \ref{lem:det}.

From either (\ref{eq:posdef2}) or (\ref{eq:condJ3}) it follows that $L_{1}-\alpha_{1}$ and $L_{2}-\alpha_{2}$
have the same sign and also that $\frac{1}{\alpha_{1}}-J_{1}$ and $\frac{1}{\alpha_{2}}-J_{2}$ have the same sign.

A straight forward calculation shows

\begin{align}
L_{1}~ >~\alpha_{1}\quad&\iff\quad
J_{1}-\frac{1}{\alpha_{1}}~<~\frac{\bar{J}^{2}}{J_{2}}\label{J1}\\
\text{and}\qquad L_{2}~ >~\alpha_{2}\quad&\iff\quad
J_{2}-\frac{1}{\alpha_{2}}~<~\frac{\bar{J}^{2}}{J_{1}}\label{J2} \, .
\end{align}

Now suppose \eqref{eq:posdef1} and \eqref{eq:posdef2} hold. If $J_{1}-\frac{1}{\alpha_{1}}>0$ holds, then also $J_{2}-\frac{1}{\alpha_{2}}> 0$.

Consequently \eqref{J1} and \eqref{J2} imply that
\begin{align}
(J_{1}-\frac{1}{\alpha_{1}})(J_{2}-\frac{1}{\alpha_{2}})~<~\frac{\bar{J}^{4}}{J_{1}J_{2}}~<~\bar{J}^{2}.
\end{align}
but this contradicts \eqref{eq:condJ3}. Thus $\frac{1}{\alpha_{1}}-J_{1}>0$.

If, on the other hand, $\frac{1}{\alpha_{1}}-J_{1}>0$ then \eqref{J1} implies $L_{1}>\alpha_{1}$.

\end{proof}
\subsubsection{Critical regime}
We turn to the critical regime:
\begin{prop}
If $L_{\nu}-\alpha_{\nu}>0$ for both groups and $(L_{1}-\alpha_{1})(L_{2}-\alpha_{2})=\bar{L}^{2}$,
then the function F defined in \eqref{eq:def_F} has a unique global
minimum at the origin.
\end{prop}

\begin{rem}
The conditions stated in the proposition are equivalent to the critical
regime. This is shown in analogous fashion to the proof of Proposition
\ref{prop:equivalence}.
\end{rem}

\begin{proof}
We take derivatives of $F$ with respect to both variables

\begin{align*}
F_{1}(y_{1},y_{2}) & =L_{1}y_{1}-\bar{L}y_{2}-\alpha_{1}\tanh y_{1}=0,\\
F_{2}(y_{1},y_{2}) & =L_{2}y_{2}-\bar{L}y_{1}-\alpha_{2}\tanh y_{2}=0.
\end{align*}

One solution to this system of equations is $y_{1}=y_{2}=0$. We proceed
to show that this solution is unique. We rewrite the function $F$:

\[
F(tx_{0},ty_{0})=\frac{1}{2}L_{1}t^{2}x_{0}^{2}+\frac{1}{2}L_{2}t^{2}y_{0}^{2}-\bar{L}x_{0}y_{0}t^{2}-\alpha_{1}\ln\cosh tx_{0}-\alpha_{2}\ln\cosh ty_{0},
\]
where $(x_{0},y_{0})$ indicates the direction, $x_{0}^{2}+y_{0}^{2}=1$,
and $t$ is the distance from the origin. The first derivative of
$F$ with respect to $t$ is 0 at the origin, independently of the
direction $(x_{0},y_{0})$ .

We show that the second derivative $\frac{\mathrm{d}^{2}F(tx_{0},ty_{0})}{\mathrm{d}t^{2}}$
is positive in all directions, except for two.

\begin{align*}
\frac{\mathrm{d}^{2}F(tx_{0},ty_{0})}{\mathrm{d}t^{2}} & =L_{1}x_{0}^{2}+L_{2}y_{0}^{2}-2\bar{L}x_{0}y_{0}-\frac{\alpha_{1}x_{0}^{2}}{\cosh^{2}tx_{0}}-\frac{\alpha_{2}y_{0}^{2}}{\cosh^{2}ty_{0}}.
\end{align*}

Therefore, we have

\begin{align*}
\left.\frac{\mathrm{d}^{2}F(tx_{0},ty_{0})}{\mathrm{d}t^{2}}\right|_{t=0} & \geq0
\end{align*}
with equality if and only if both $t=0$ and

\[
\sqrt{L_{1}-\alpha_{1}}x_{0}-\sqrt{L_{2}-\alpha_{2}}y_{0}=0
\]
hold.

Hence there are two directions $(x_{0},y_{0})$, one pointing into
quadrant one, the other into quadrant three, in which the second derivative
is 0 at the origin. In all other directions the second derivative
is strictly positive. For any direction, the second derivative is
strictly positive for all $t>0$.

This concludes the proof that the minimum at the origin is unique
and global.
\end{proof}
\subsubsection{Low Temperature Regime} 
Assume we are in the low temperature regime, i.e. at least one of
the following conditions holds:

\begin{align}
J_{1}~ & >~\frac{1}{\alpha_{1}},\label{eq:condJ1-2}\\
J_{2}~ & >~\frac{1}{\alpha_{2}},\label{eq:condJ2-2}\\
\bar{J}^{2}~ & >~(\frac{1}{\alpha_{1}}-J_{1})(\frac{1}{\alpha_{2}}-J_{2}).\label{eq:condJ3-2}
\end{align}

In terms of the inverse matrix $L=J^{-1}$, at least one of the following
inequalities has to hold:

\begin{align}
L_{1} & \leq\alpha_{1},\label{eq:condL1}\\
L_{2} & \leq\alpha_{2},\label{eq:condL2}\\
(L_{1}-\alpha_{1})(L_{2}-\alpha_{2}) & <\bar{L}^{2}.\label{eq:condL3}
\end{align}

In order to apply Laplace's method, we need to determine the minima
of the function

\begin{equation}
F(x,y)=\frac{1}{2}L_{1}x^{2}+\frac{1}{2}L_{2}y^{2}-\bar{L}xy-\alpha_{1}\ln\cosh x-\alpha_{2}\ln\cosh y.
\end{equation}

The first order conditions are
\begin{align}
F_{1}(x,y) & =L_{1}x-\bar{L}y-\alpha_{1}\tanh x=0,\label{eq:FOC1}\\
F_{2}(x,y) & =L_{2}y-\bar{L}x-\alpha_{2}\tanh y=0.\label{eq:FOC2}
\end{align}

These equations always have a solution $(x,y)=(0,0)$. We define functions
$X:[0,\infty)\rightarrow[0,\infty)$ and $Y:[0,\infty)\rightarrow[0,\infty)$
by setting $X(y)$ equal to the largest solution $x$ of equation
\eqref{eq:FOC1} given a value $y\geq0$. Similarly, $Y(x)$ is defined
as the largest solution $y$ of \eqref{eq:FOC2} given $x\geq0$.
\begin{prop}
The functions $X$ and $Y$ are strictly increasing and strictly concave.
\end{prop}

\begin{proof}
We show the properties for $X$. By the implicit function theorem,
we can calculate the first derivative of the function $X$ by dividing
the partial derivative of the function

\[
G(x,y):=L_{1}x-\bar{L}y-\alpha_{1}\tanh x
\]

with respect to $x$ by the partial derivative of $G$ with respect
to $y$. That yields

\[
X'(y)=\frac{\bar{L}}{L_{1}-\frac{\alpha_{1}}{\cosh^{2}X(y)}}.
\]

We show that $X'$ is always positive. We define two auxiliary functions

\begin{align*}
f,g & :[0,\infty)\rightarrow[0,\infty),\\
f(x) & :=L_{1}x-\bar{L}y,\\
g(x) & :=\alpha_{1}\tanh x.
\end{align*}

We are looking for the intersections of the functions $f$ and $g$
given a value of $y\geq0$, the largest of which is precisely the
value $X(y)$. If $y=0$ and $L_{1}-\alpha_{1}\geq0$, $f$ and $g$
only intersect at 0, so $X(0)=0$ in this case. If $L_{1}-\alpha_{1}<0$,
or $y>0$ holds, then $f(0)=-\bar{L}y\leq0$, $g(0)=0$, and $f'(0)<g'(0)$,
so at the origin $f(0)\leq g(0)$, and for small values of $x$ $f(x)<g(x)$.
However, whereas $f'(x)=L_{1}>0$ is constant, $g'(x)=\frac{\alpha_{1}}{\cosh^{2}x}>0$
is strictly decreasing in $x$ and $\lim_{x\rightarrow\infty}g'(x)=0$.
Therefore there is exactly one value $x_{1}>0$ such that $f(x_{1})=g(x_{1})$.
Since for $x<x_{1}$ $f(x)<g(x)$, it must be that $f'(x_{1})>g'(x_{1})$.
This $x_{1}$ is $X(y)$. Hence we have

\[
L_{1}>\frac{\alpha_{1}}{\cosh^{2}X(y)},
\]

and $X'(y)>0$ has been shown.

The second derivative of $X$ is

\[
X''(y)=-\frac{2\alpha_{1}\bar{L}X'(y)}{\left(L_{1}-\frac{\alpha_{1}}{\cosh^{2}X(y)}\right)^{2}\cosh^{3}X(y)}<0,
\]

and so $X$ is strictly concave.
\end{proof}
\begin{prop}
The limits of the first derivatives of $X$ and $Y$ are

\begin{align*}
\lim_{y\rightarrow\infty}X'(y) & =\frac{\bar{L}}{L_{1}}>0,\\
\lim_{x\rightarrow\infty}Y'(x) & =\frac{\bar{L}}{L_{2}}>0.
\end{align*}
\end{prop}

\begin{proof}
We have

\begin{align*}
\lim_{y\rightarrow\infty}X'(y) & =\frac{\bar{L}}{L_{1}-\frac{\alpha_{1}}{\cosh^{2}X(y)}}\\
 & =\frac{\bar{L}}{L_{1}}
\end{align*}

due to $\frac{\alpha_{1}}{\cosh^{2}X(y)}\rightarrow0$ as $X(y)\rightarrow\infty$.
On the other hand, as $y$ goes to infinity, the solution $X(y)$
of \eqref{eq:FOC1} has to go to infinity due to the boundedness of
the term $\alpha_{1}\tanh x$.
\end{proof}
\begin{cor}
\begin{align}
\lim_{y\rightarrow\infty}X'(y) & \lim_{x\rightarrow\infty}Y'(x)=\frac{\bar{L}^{2}}{L_{1}L_{2}}<1\label{eq:prod_lim_derivatives}
\end{align}
 holds.
\end{cor}

\begin{proof}
By assumption, $J$, and therefore $L=J^{-1}$, are positive definite.
In particular, the determinant of $L$ must be positive. Hence

\[
L_{1}L_{2}-\bar{L}^{2}>0.
\]
\end{proof}
We define the curves $\gamma_{1},\gamma_{2}:[0,\infty)\rightarrow\mathbb{R}^{2}$
by setting

\begin{align*}
\gamma_{1}(y) & :=(X(y),y),\\
\gamma_{2}(x) & :=(x,Y(x)).
\end{align*}

These curves originate at a certain point in the first quadrant that
depends on the parameters $L_{1},L_{2},\alpha_{1},\alpha_{2}$:

\begin{align*}
\gamma_{1}(0) & :=\begin{cases}
(0,0), & L_{1}-\alpha_{1}\geq0,\\
(X(0),0), & L_{1}-\alpha_{1}<0,
\end{cases}\\
\gamma_{2}(x) & :=\begin{cases}
(0,0), & L_{2}-\alpha_{2}\geq0,\\
(0,Y(0)), & L_{2}-\alpha_{2}<0.
\end{cases}
\end{align*}

So curve $\gamma_{1}$ starts at the origin if and only if $L_{1}-\alpha_{1}\geq0$.
Similarly, $\gamma_{2}$ starts at the origin if and only if $L_{1}-\alpha_{1}\geq0$.

Let us first assume the two curves do not meet at the origin. Then
they start at points apart, but due to \eqref{eq:prod_lim_derivatives},
we have

\[
\lim_{x\rightarrow\infty}Y'(x)<\frac{1}{\lim_{y\rightarrow\infty}X'(y)}.
\]

This implies the curves have to meet at some point in the interior
of the first quadrant. Call this point $(x_{1},y_{1})$. Once they
have met, the strict concavity of both $X$ and $Y$ drives them apart
and they do not intersect again. Hence the point $(x_{1},y_{1})$
is uniquely determined.

If the two curves do meet at the origin, we have $L_{1}-\alpha_{1}\geq0$
and $L_{2}-\alpha_{2}\geq0$. We distinguish the two cases
\begin{enumerate}
\item $L_{1}-\alpha_{1}>0$ and $L_{2}-\alpha_{2}>0$,
\item $L_{1}-\alpha_{1}=0$ or $L_{2}-\alpha_{2}=0$.
\end{enumerate}
In the first case, since we are in the low temperature regime,

\[
(L_{1}-\alpha_{1})(L_{2}-\alpha_{2})<\bar{L}^{2}
\]

must hold, and we have
\begin{align*}
X'(0)Y'(0) & =\frac{\bar{L}}{L_{1}-\alpha_{1}}\frac{\bar{L}}{L_{2}-\alpha_{2}}\\
 & >1.
\end{align*}

So

\[
Y'(0)>\frac{1}{X'(0)},
\]

which means that the two curves starting at $(0,0)$ initially move
apart. So for some points close to the origin, the curves are apart,
and the previous reasoning for the existence of a unique point of
intersection inside the first quadrant applies.

In the second case, if $L_{1}-\alpha_{1}=0$, then the derivative
of function $X$ is infinite at the origin, meaning $\gamma_{1}$
moves parallel to the x-axis. The function $Y$ on the other hand
either has positive or infinite derivative at the origin, so $\gamma_{2}$
either moves into the interior of the first quadrant or it moves parallel
to the y-axis. In any case, the two curves move apart after leaving
the origin.

We summarize:

\begin{thm}\label{thm:corrl}
   In the low temperature regime the function $F$ has exactly two minima $\mu^{*}=(\mu_{1}^{*}, \mu_{2}^{*})$ and $-\mu^{*}$
   and we may suppose that $\mu_{1}^{*}, \mu_{2}^{*}>0$.
\end{thm}

\subsection{Correlations for the high temperature regime}\label{sub:corrhigh}

In this section we use the Laplace method to evaluate the expression $Z_{J}(K,Q)$ and thus the correlation \eqref{eq:correlations} asymptotically in the high temperature regime.

Let $H=J^{-1}- \begin{pmatrix}\alpha_{1}&0\\0&\alpha_{2}\end{pmatrix}$
be the Hessian of $F$ at $0$, $\mathcal{N}(0,H^{-1})$ the two-dimensional normal distribution with covariance matrix $H^{-1}$
and let
\begin{align}
\nu_{K,Q}~&=~\nu_{K,Q}(0,H^{-1})\notag\\~&=~\frac{\sqrt{\det H}}{2\pi}\,\int_{\IR^{2}}e^{-\frac{1}{2}
\begin{pmatrix}x_{1}& x_{2}\end{pmatrix}H\begin{pmatrix}x_{1}\\x_{2}\end{pmatrix}} {x_{1}}^{K}\,{x_{2}}^{Q}\,dx_{1}\,dx_{2}
\end{align}
be the moments of $\mathcal{N}(0,H^{-1})$.

In the following proposition as in the whole paper by $a_{N}\approx b_{N}$ we mean $\lim_{N\to\infty} \frac{a_{N}}{b_{N}}=1$.

\begin{prop}\label{prop:ZJhigh}
   Let $J,\alpha_{1},\alpha_{2}$ satisfy \eqref{eq:condJ1-1}--\eqref{eq:condJ3-1} (high temperature regime), then
   \begin{align}
       Z_{J}(K,Q)~&=~\int e^{-N F_{J}(y)}\;\tanh^{K}y_{1}\,\tanh^{Q}y_{2}\,\mathrm{d}^{2}y\notag\\
       &\approx~\frac{2\pi}{\sqrt{\det H}}\;\nu_{K,Q}(0,H^{-1})\;N^{-\frac{K+Q}{2}-1}\qquad \text{as $N\to\infty$}\,.
   \end{align}
\end{prop}

\begin{proof}
We use the Laplace method to evaluate $Z_{J}(K,Q)$. We only sketch the main idea. For the details, in particular the remainder estimates, we refer to \cite{Olver} or \cite{MM}. In the integral we replace both $F$ and the tanh terms by the leading terms in their Taylor expansion around 0.
This gives
 \begin{align*}
  Z_{J}(K,Q)
 & \approx \int e^{-N \cdot 1/2(y'Hy)} \,  y_1^K \, y_2^Q \, dy_1 \, dy_2  \\
 & = N^{-\frac{K+Q}{2}-1} \int e^{-\frac{1}{2}(x'Hx)} \, x_1^K \, x_2^Q \, dx_1 \, dx_2 \\
      & = \frac{2\pi}{\sqrt{\det H}}\;\nu_{K,Q}(H^{-1})\;N^{-\frac{K+Q}{2}-1}
   \end{align*}
   where we changed variable $x= \sqrt{N}y$.
\end{proof}

Proposition \ref{prop:ZJhigh} immediately gives:
\begin{thm}\label{thm:mohi}
Let $J, \alpha_1, \alpha_2$ be in the high temperature regime (\eqref{eq:condJ1-1}--\eqref{eq:condJ3-1} then
\begin{equation}
\IE(X_1 \cdot X_2 \cdot ... \cdot X_K \cdot Y_1 \cdot Y_2 \cdot ... \cdot Y_Q) \approx \nu_{K,Q}(0, H^{-1})\; N^{-(K+Q)/2}
\end{equation}
\end{thm}

\subsection{Correlations for the critical regime}\label{sub:corrcrit}

Expanding again $F$ to leading order gives in the critical regime

\begin{align*}
  F(y_1,y_2) &\approx \frac{1}{2} ((L_1 - \alpha_1) y_1^2 + (L_2 - \alpha_2)y_2^2 - 2  \bar{L}y_1y_2 + \frac{2 \alpha_1}{12}y_1^4+ \frac{2 \alpha_2}{12}y_2^4) \\
 & =  \frac{1}{2} ( \, (\sqrt{L_1 - \alpha_1} y_1 - \sqrt{L_2- \alpha_2}y_2)^2 + \frac{\alpha_1}{6} y_1^4 + \frac{\alpha_2}{6}y_2^4 \,)
 \end{align*}

Thus we have

\begin{align*}
  Z_{J}(K,Q)  \approx \int_{\IR^2} e^{-N/2( \, (\sqrt{L_1 - \alpha_1} y_1 - \sqrt{L_2- \alpha_2}y_2)^2 + \frac{\alpha_1}{6} y_1^4 + \frac{\alpha_2}{6}y_2^4 \,) } y_1^K \, y_2^Q \, dy_1 \, dy_2
\end{align*}

 We substitute

 \begin{align*}
 u & = N^{1/2} (\sqrt{L_1 - \alpha_1} \,y_1 - \sqrt{L_2- \alpha_2} \, y_2) \\
 v & = N^{1/4} (\sqrt{L_1 - \alpha_1} \,y_1 + \sqrt{L_2- \alpha_2} \, y_2)
 \end{align*}

 which gives
 \[
\int_{\mathbb{R}^{2}}e^{-\frac{1}{2}\left[u'^{2}+\frac{\alpha_{1}}{2^{5}\cdot3(L_{1}-\alpha_{1})^{2}}(\frac{u'}{N^{1/4}}+v')^{4}+\frac{\alpha_{2}}{2^{5}\cdot3(L_{2}-\alpha_{2})^{2}}(v'-\frac{u'}{N^{1/4}})^{4}\right]}(\frac{u'}{N^{1/2}}+\frac{v'}{N^{1/4}})^{K}(\frac{v'}{N^{1/4}}-\frac{u'}{N^{1/2}})^{L}\mathrm{d}u'\mathrm{d}v'
\]
times a constant equal to

\[
\frac{1}{2^{K+L+1}(L_{1}-\alpha_{1})^{\frac{K+1}{2}}(L_{2}-\alpha_{2})^{\frac{L+1}{2}}N^{\frac{3}{4}}}.
\]

Since we are interested merely in the \emph{ratio} $\frac{Z_{J}(K,Q)}{Z_{J}(0,0)}$ we may (and will) neglect multiplicative constants in the evaluation of $Z_{J}(K,Q) $ as long as
these constants are independent of $K$ and $Q$. To shorten notation we define
\begin{align}
   a_{N}(K,Q)~\sim~b_{N}(K,Q) \qquad\text{ if }\qquad \lim_{N\to\infty}\frac{a_{N}}{b_{N}}~\to~c
\end{align}
for a constant $0<c<\infty$ which is independent of $K$ and $Q$. With this notation we have

\begin{align*}
  Z_{J}(K,Q) & \sim \int_{\IR^2} e^{-1/2u^2} \, e^{-1/(3 \cdot 26) \left( \frac{\alpha_1}{(L_1 - \alpha_1)^2} +  \frac{\alpha_2}{(L_2-\alpha_2)^2)} v^4 \right) } v^{K+Q} \, du \, dv \\
 & \sim     \int e^{-1/(3 \cdot 26) \left( \frac{\alpha_1}{(L_1 - \alpha_1)^2} +  \frac{\alpha_2}{(L_2-\alpha_2)^2)} v^4 \right) } v^{K+Q} \, dv
\end{align*}
for constants $c, \, c'$ independent of $K$ and $Q$. We note that
\begin{equation*}
\int_0^\infty e^{-a\cdot x^4} \, x^m \, dx = \frac{1}{4a^\frac{m+1}{4}} \Gamma(\frac{m+1}{4})
\end{equation*}
Summing up, we obtain
\begin{thm}\label{thm:mocri}
   Let $J, \alpha_{1}, \alpha_{2}$ be in the critical regime (\eqref{eq:condJ1-1-1} -- \eqref{eq:condJ2-1-1}) then
\begin{align}
&\IE ( X_1 \cdot ... \cdot X_K \cdot Y_1 \cdot ... \cdot Y_Q)\notag
 \approx \left[\frac{12}{\alpha_1 (L_2 -\alpha_2)^2 + \alpha_2(L_1 -\alpha_1)^2} \right]^\frac{K+Q}{4}\notag
\\&\hspace*{12pt}
\cdot (L_1 - \alpha_1)^{Q/2} (L_2 - \alpha_2)^{K/2} \frac{\Gamma(\frac{K+Q+1}{4})}{\Gamma(\frac{1}{4})}  \cdot N^{-\frac{K+Q}{4}}\label{eq:mocri}
\end{align}
\end{thm}

\subsection{Correlations for the Low Temperature Regime}
Using again Laplace's method to evaluate the expressions for $Z_{J}(K,Q)$ for the low temperature regime we immediately get:

\begin{thm}
   In the low temperature regime we have
   \begin{align}
      \IE\Big(X_1 \cdot ... \cdot X_K \cdot Y_1 \cdot ... \cdot Y_Q)\Big)~\approx~\left\{
                                                                                    \begin{array}{ll}
                                                                                      \tanh^{K}(\mu_{1})\,\tanh^{Q}(\mu_{2}), & \hbox{if $K+Q $ is even;} \\
                                                                                      0, & \hbox{otherwise.}
                                                                                    \end{array}
                                                                                  \right.
             \end{align}
    where $\mu_{*}=(\mu^{*}_{1},\mu^{*}_{2})$ is given in Theorem \ref{thm:corrl}.
\end{thm}

\section{Laws of large numbers}\label{sec:lln}
In this section we prove Theorems \ref{LLN}, \ref{LLN-1}, and \ref{thm:SpecialCase}.
We set
\begin{align}
   W_{K,N_{1}}~:=~\{1,2,\ldots,N_{1}\}^{K}\,.
\end{align}
We also denote by $W_{K,N_{1}}(r)$ the set of all multiindices $\ui=(i_{1},i_{2},\ldots,i_{K})\in W_{K,N_{1}}$
for which exactly $r$ indices occur only once and by $w_{K,N_{1}}(r)$ the cardinality of $W_{K,N_{1}}(r)$.

We have
\begin{lem}\label{thm:comb1}
\begin{align}
   w_{K,N_{1}}(r)~\leq~K!\,N_{1}^{\frac{K+r}{2}}
\end{align}
\end{lem}

\begin{proof}
   The multiindices in $W_{K,N_{1}}(r)$ contain at most $r+\frac{K-r}{2}=\frac{K+r}{2}$ different indices. There are at most $N_{1}^{\frac{K+r}{2}} $ ways to choose them
   and at most $K!$ ways to order them.
\end{proof}

\begin{thm}\label{thm:llnmoments}
   If \eqref{eq:condJ1-1} --\eqref{eq:condJ3-1} (high temperature regime) hold, then for all $K,Q\in\IN$, $K,Q>0$
   \begin{align}\label{eq:lln-moment}
\IE\Big(\big(\frac{1}{N_{1}}\sum_{i=1}^{N_{1}}X_{i}\big)^{K}\;\big(\frac{1}{N_{2}}\sum_{j=1}^{N_{2}}Y_{j}\big)^{Q}\Big)~\rightarrow~0
   \end{align}
   \eqref{eq:lln-moment} is also true if \eqref{eq:condJ1-1-1} --\eqref{eq:condJ3-1-1} (critical regime) hold.
\end{thm}
\begin{proof}
   \begin{align}
      &\IE\Big(\big(\frac{1}{N_{1}}\sum_{i=1}^{N_{1}}X_{i}\big)^{K}\;\big(\frac{1}{N_{2}}\sum_{j=1}^{N_{2}}Y_{j}\big)^{Q}\Big)\notag\\
      =~&\frac{1}{N_{1}^{K}\,N_{2}^{Q}} \sum_{\ui\in W_{K,N_{1}}}\;\sum_{\uj\in W_{Q,N_{2}}}\,
      \IE\Big(X_{i_{1}}\cdot X_{i_{2}}\cdot\ldots\cdot X_{i_{K}}\,Y_{j_{1}}\cdot Y_{j_{2}}\cdot\ldots\cdot Y_{j_{Q}}\Big)\notag\\
      =~&\frac{1}{N_{1}^{K}\,N_{2}^{Q}} \sum_{k=0}^{K}\,\sum_{q=0}^{Q}\;\sum_{\ui\in W_{K,N_{1}}(k)}\;\sum_{\uj\in W_{Q,N_{2}}(q)}\,
      \IE\Big(X_{i_{1}}\cdot X_{i_{2}}\cdot\ldots\cdot X_{i_{K}}\,Y_{j_{1}}\cdot Y_{j_{2}}\cdot\ldots\cdot Y_{j_{Q}}\Big)\notag\\
      \leq~&C\;\frac{1}{N_{1}^{K}\,N_{2}^{Q}}\,N_{1}^{\frac{K+k}{2}}\,N_{2}^{\frac{Q+q}{2}}\,N^{-(k+q)/2}~\to 0\label{eq:llnest}
   \end{align}
   where we used Theorem \ref{thm:mohi} and \eqref{eq:lln-moment} in the final estimate. Note, that the constant $C$ depends on $K$ and $Q$, but not on $N$.

   Using the estimate \eqref{eq:mocri} instead we obtain the result in the critical regime as well.

      We also note that the estimates hold in the cases $\alpha_{1}=0$ or $\alpha_{2}=0$.
\end{proof}
Theorem \ref{LLN} and Theorem \ref{LLN-1} follow immediately from \ref{thm:llnmoments} and \ref{thm:moments}.
\begin{rem}
   We have actually proved that
   \begin{align*}
      \IE\Big(\big(\frac{1}{N_{1}^{\gamma}}\sum_{i=1}^{N_{1}}X_{i}\big)^{K}\;\big(\frac{1}{N_{2}^{\gamma}}\sum_{j=1}^{N_{2}}Y_{j}\big)^{Q}\Big)
   \end{align*}
   is bounded with $\gamma=\frac{1}{2}$ in the high temperature regime and $\gamma=\frac{3}{4}$ in the critical regime. This is an indication that
   the limit theorems \ref{CLT} and \ref{Fluctuations} may hold.
\end{rem}

We turn to the low temperature regime. Analogous to \eqref{eq:llnest} we obtain
 \begin{align}
      &\IE\Big(\big(\frac{1}{N_{1}}\sum_{i=1}^{N_{1}}X_{i}\big)^{K}\;\big(\frac{1}{N_{2}}\sum_{j=1}^{N_{2}}Y_{j}\big)^{Q}\Big)\notag\\
           =~&\frac{1}{N_{1}^{K}\,N_{2}^{Q}} \sum_{k=0}^{K}\,\sum_{q=0}^{Q}\;\sum_{\ui\in W_{K,N_{1}}(k)}\;\sum_{\uj\in W_{Q,N_{2}}(q)}\,
      \IE\Big(X_{i_{1}}\cdot X_{i_{2}}\cdot\ldots\cdot X_{i_{K}}\,Y_{j_{1}}\cdot Y_{j_{2}}\cdot\ldots\cdot Y_{j_{Q}}\Big)\label{eq:llnlow}
   \end{align}
The terms with $k<K$ and $q<Q$ are canceled by the term in front of the sum. Thus \eqref{eq:llnlow} is asymptotically given by
\begin{align}
   &\IE\Big(X_{1}\cdot X_{2}\cdot\ldots\cdot X_{K}\,Y_{1}\cdot Y_{2}\cdot\ldots\cdot Y_{Q}\Big)\notag\\
~\approx~&\frac{1}{2}\big(1+(-1)^{K+Q}\big)\; \tanh^{K}(\mu^{*}_{1})\cdot\tanh^{Q}(\mu^{*}_{2})
\end{align}

\section{The Central Limit Theorem}
\subsection{Some Combinatorics}\label{subs:comb}

To prove the Central Limit Theorem \ref{CLT} we need a more detailed analysis of Lemma \ref{thm:comb1}.

Let us define $W_{K,N_{1}}^{0}(r)$ to be the set of all multiindices $\ui=(i_{1},i_{2},\ldots,i_{K})\in W_{K,N_{1}}(r)$
for which no index occurs more than twice. We also set
$$W_{K,N_{1}}^{+}(r):=W_{K,N_{1}}(r)\setminus W_{K,N_{1}}^{0}(r)$$
and denote by $w_{K,N_{1}}^{+}(r)$ and $w_{K,N_{1}}^{0}(r)$ their cardinalities.

\begin{lem}\label{lem:comb2}
\begin{align}
   w_{K,N_{1}}^{+}(r)~\leq~K!\,N_{1}^{\frac{K+r}{2}-\frac{1}{2}}\,.
\end{align}
\end{lem}

\begin{proof}
   If the $K$-tuple $\ui$ contains $r$ indices with only one occurrence and at least one index with three or more occurrences
   there are at most $r-3$ places left for indices with (exactly) two occurrences. Therefore, a tuple in $w_{K,N_{1}}^{+}(r)$
   contains at most $k + 1 +\frac{K-r-3}{2}$ \emph{different} indices. Consequently there are at most
   $K!\,N_{1}^{\frac{K+r}{2}-\frac{1}{2}}$ such tuples.
\end{proof}

\begin{lem}\label{lem:comb4}
\begin{align}
   w_{K,N_{1}}^{0}(r)~=~\left\{
                          \begin{array}{ll}
                           \frac{N_{1}!}{(N_{1}-\frac{K+r}{2})!}\;\frac{K!}{r!\;(\frac{K-r}{2})!\; 2^{\frac{K-r}{2}}}\; , & \hbox{if $K-r$ is even;} \\
                            0, & \hbox{else.}
                          \end{array}
                        \right.
  \end{align}
\end{lem}
\begin{proof}
   We choose an (ordered) $r$-tuple $\rho$ of $r$ indices to occur once and an ordered  $(K-r)/2$-tuple $\lambda$ of indices to occur twice in $\ui$.
We have
\begin{align*}
   \frac{N_{1}!}{(N_{1}-\frac{K+r}{2})!}
\end{align*}
ways to do so.

Then we choose the $r$ positions for those indices which occur once. We can do this in
\begin{align*}
   \binom{K}{r}~=~\frac{K!}{r!\;(K-r)!}
\end{align*}
ways. We fill these positions in $\ui$ with $\rho_{1}, \rho_{2},\ldots,\rho_{r}$ starting with the left most open position.

Finally we distribute the indices $\lambda_{1},\ldots,\lambda_{(K-r)/2}$, twice each. The index $\lambda_{1}$ is put at the left most free place in $\ui$ and in
one of the remaining $K-r-1$ positions, $\lambda_{2}$ is put at the then first free place in $\ui $ and in one of the $K-r-3$ remaining free places and so on.

This gives
\begin{align}
   (K-r-1)!!~=~\frac{(K-r)!}{(\frac{K-r}{2})!\; 2^{\frac{K-r}{2}}}
\end{align}
 possibilities.
\end{proof}

We summarize the above considerations in the following Theorem.

\begin{thm}\label{thm:comb3}
In the high temperature regime we have

\begin{align}
   &\IE\Big(\big(\frac{1}{N_{1}^{\frac{1}{2}}}\sum_{i=1}^{N_{1}}X_{i}\big)^{2K}\;\big(\frac{1}{N_{2}^{\frac{1}{2}}}\sum_{j=1}^{N_{2}}Y_{j}\big)^{2Q}\Big)\notag\\
   \approx~&\sum_{k=0}^{K}\,\sum_{q=0}^{Q}\;\frac{(2K)!\alpha_1^k}{(2k)!\;(K-k)!\; 2^{K-k}}\;\frac{(2Q)!\alpha_2^q}{(2q)!\;(Q-q)!\; 2^{Q-q}}\ \nu_{2k,2q}(0,H^{-1})\label{eq:comb3-1}\\
   \intertext{and}
   &\IE\Big(\big(\frac{1}{N_{1}^{\frac{1}{2}}}\sum_{i=1}^{N_{1}}X_{i}\big)^{2K+1}\;\big(\frac{1}{N_{2}^{\frac{1}{2}}}\sum_{j=1}^{N_{2}}Y_{j}\big)^{2Q+1}\Big)\notag\\
   \approx~&\sum_{k=0}^{K}\,\sum_{q=0}^{Q}\;\frac{(2K+1)!\alpha_1^{k+\frac{1}{2}}}{(2k+1)!\;(K-k)!\; 2^{K-k}}\;\frac{(2Q+1)!\alpha_2^{q+\frac{1}{2}}}{(2q+1)!\;(Q-q)!\; 2^{Q-q}}\ \nu_{2k+1,2q+1}(0,H^{-1})\label{eq:comb3-2}
\end{align}
\end{thm}

\subsection{Moments of a 2d-normal distribution}
Let $\Sigma =  \left(\begin{array}{cc} \sigma_1 & \bar{\sigma}\\ \bar{\sigma} & \sigma_2
\end{array}\right) $ be a covariance matrix and let  $ \left(\begin{array}{c} Z_1 \\ Z_2
\end{array}\right) \, \sim \mathcal{N}(0,\Sigma)$ distributed.

We write $ \nu_{K,Q} (\Sigma) = \mathbb{E}(Z_1^K Z_2^Q)$ to denote the moment of $\mathcal{N}(0, \Sigma)$ of order $(K,Q)$.

\begin{prop}\label{prop:Isserlis}
\begin{equation*}
\nu_{2K,2Q}(\Sigma)= \sum_{r=0}^{K \wedge Q} \frac{2K!}{(2K - 2r)! \, 2r!} \, \frac{2Q!}{(2Q - 2r)!} \frac{(2K - 2r)!}{(K-r)! \, 2^{K-r}} \frac{(2Q-2r)!}{(Q-r)! \, 2^{Q-r}} \sigma_1^{K-r} \bar{\sigma}^{2r} \sigma_2^{Q-r}
\end{equation*}
\end{prop}

\begin{proof}
Suppose $V_{1},\ldots,V_{2n}$ are random variables. Denote by $\mathcal{P}_{2}=\mathcal{P}_{2}(2n)$ the set of pair partitions of $1,\ldots, 2n$.
For a pair partition $\pi=\{\pi_{1},\ldots,\pi_{n}\}\in\mathcal{P}_{2}$ we set $\prod_{\pi}(V_{1},\ldots,V_{2n})=\prod_{i=1}^{n}\IE(V^{\pi_{i}})$
where $V^{\{i,j\}}=V_{i}\cdot V_{j}$.

By Isserlis' theorem \cite{Isserlis} we have \\
\begin{align*}
\nu_{2K,2Q}(\Sigma) = & \sum_{\pi \in \mathcal{P}_2(2K+2Q} \Pi_{\pi}\Big(\underbrace{Z_{1},Z_{1},\ldots,Z_{1}}_{2K \text{ times}}\,\underbrace{Z_{2},Z_{2},\ldots,Z_{2}}_{2Q \text{ times}}\Big) \\
= & \sum_{r=0}^{K \wedge Q} \rho_{r} \cdot \sigma_1^{2K - 2r} \bar{\sigma}^{2r} \sigma_2^{2Q - 2r}
\end{align*}
where $\rho_{r}=\# \lbrace \pi | \pi \text{ contains exactly }2r\text{ mixed pairs} \rbrace$. (Mixed pairs are of the form $\{i,j\}$ with $i\leq 2K$ and $j> 2K$.)

To compute $\rho_{r} $ we first choose $2r$ "$Z_1$'s". This can be done in $\binom{2K}{2r}$ ways.
For these  $Z_1$'s choose $2r$  $Z_2$'s :  $ \frac{(2L)!}{(2L- 2r)!}$. This gives $\left(\begin{array}{c} 2K \\ 2r \end{array}\right) \frac{(2L)!}{(2L- 2r)!}$.
The remaining terms come from the pair partitions of the $2K-2r$ $Z_1$'s and $2Q-2R$ $Z_2$'s.

\end{proof}

\subsection{Proof of Theorem \ref{CLT}}\label{sub:prclt}
We calculate the moments $ \mathbb{E} \left( (\frac{1}{N_1^{1/2}} \sum X_i)^{2K} \, (\frac{1}{N_2^{1/2}} \sum Y_j)^{2Q} \right) $, i.~e. those with even exponents. The case of
odd exponents is done in a similar way.

We write $H^{-1}=\begin{pmatrix}
                    \sigma_{1} & \bar{\sigma}\\ \bar{\sigma} & \sigma_{2}
                 \end{pmatrix}$

From Theorem \ref{thm:comb3} and Proposition \ref{prop:Isserlis} we know that the moments  $ \mathbb{E} \left( (\frac{1}{N_1^{1/2}} \sum X_i)^{2K} \, (\frac{1}{N_2^{1/2}} \sum Y_j)^{2Q} \right) $ are approximately given by

\begin{align*}
& \sum_{k=0}^K \sum_{l=0}^Q \frac{2K!\alpha_1^k}{2k! \, (K-k)! \, 2^{K-k}} \, \frac{2Q! \alpha_2^l}{2l! \, (Q-l)! \, 2^{Q-l}} \, \nu_{k,l} (H^{-1}) \\
=& \sum_{k=0}^K \sum_{l=0}^Q \frac{2K!\alpha_1^k}{2k! \, (K-k)! \, 2^{K-k}} \, \frac{2Q! \alpha_2^l}{2l! \, (Q-l)! \, 2^{Q-l}} \\ & \cdot
\sum_{r=0}^{k \wedge l} \frac{2k!}{(2k - 2r)! \, 2r!} \, \frac{2l!}{(2l - 2r)!} \frac{(2k - 2r)!}{(k-r)! \, 2^{k-r}} \frac{(2l-2r)!}{(l-r)! \, 2^{l-r}} \sigma_1^{k-r} \bar{\sigma}^{2r} \sigma_2^{l-r} \\
=& \sum_{r=0}^{K \wedge Q} \sum_{k=r}^K \sum_{l=r}^Q \frac{2K!\alpha_1^k}{(K-k)! \, 2^{K-k}} \, \frac{2Q! \alpha_2^l}{(Q-l)! \, 2^{Q-l}} \frac{2^{2r}}{2r! \, (k-r)! \, (l-r)! \, 2^k \, 2^l}  \sigma_1^{k-r} \bar{\sigma}^{2r} \sigma_2^{l-r} \\
=& \sum_{r=0}^{K \wedge Q} \sum_{k=r}^K \sum_{l=r}^Q \frac{2K!\alpha_1^k}{(K-k)! \, 2^{K-r}} \, \frac{2Q! \alpha_2^l}{(Q-l)! \, 2^{Q-r}} \frac{1}{2r! \, (k-r)! \, (l-r)! }  \sigma_1^{k-r} \bar{\sigma}^{2r} \sigma_2^{l-r}  \\
& \text{Setting }s = k-r , t = l-r \text{ (i.e. }k=s+r, l= t+r) \text{ gives:} \\
=& \sum_{r=0}^{K \wedge Q} \sum_{s=0}^{K-r} \sum_{t=0}^{Q-r} \frac{2K!\alpha_1^{s+r}}{(K-r-s)! \, 2^{K-r}} \, \frac{2Q! \alpha_2^{t+r}}{(Q-r-s)! \, 2^{Q-r}} \frac{1}{2r! \, s! \, t! } \, \sigma_1^s \bar{\sigma}^{2r} \sigma_2^t \\
=&  \sum_{r=0}^{K \wedge l} \frac{2K!}{(2K-2r)! 2r!} \frac{2Q!}{(2Q-2r)!} \frac{(2K-2r)!}{(K-r)! \, 2^{K-r}} \frac{(2Q-2r)!}{(Q-r)! 2^{Q-r}} \, ( \sqrt{\alpha_1 \alpha_2} \bar{\sigma})^{2r}        \\
& \cdot \sum_{s=0}^{K - r} \frac{(K-r)!}{(K-r-s)! \, s!} (\alpha_1 \sigma_1)^s \cdot \sum_{t=0}^{Q - r} \frac{(Q-r)!}{(Q-r-s)! \, t!} ( \alpha_2 \sigma_2)^t \\
=& 	\sum_{r=0}^{K \wedge Q} \frac{2K!}{(2K - 2r)! \, 2r!} \, \frac{2Q!}{(2Q - 2r)!} \frac{(2K - 2r)!}{(K-r)! \, 2^{K-r}} \frac{(2Q-2r)!}{(Q-r)! \, 2^{Q-r}} (1 + \alpha_1 \sigma_1)^{K-r} \,  ( \sqrt{\alpha_1 \alpha_2} \bar{\sigma})^{2r} \, (1 + \alpha_{2}\sigma_2)^{Q-r}           \\
=& \, \nu_{2K,2Q}(C)
\end{align*}
with $C$ is defined in \eqref{eq:covariance}.

This finishes the proof of Theorem \ref{CLT}.

\subsection{Proof of Theorem \ref{Fluctuations}}

Similar to the high temperature regime we evaluate
\begin{align}
    & \IE\Big(\big(\frac{1}{N_{1}^{\frac{3}{4}}}\sum_{i=1}^{N_{1}}X_{i}\big)^{K}\;\big(\frac{1}{N_{2}^{\frac{3}{4}}}\sum_{j=1}^{N_{2}}Y_{j}\big)^{Q}\Big)\notag\\
    =~&\frac{1}{N_{1}^{\frac{3}{4}\,K}\,N_{2}^{\frac{3}{4}\,Q}}\,\sum_{k=0}^{K}\sum_{q=0}^{Q} w_{K,N_{1}}^{0}(k)\; w_{Q,N_{2}}^{0}(q)\ \IE\Big(X_{1}\ldots X_{k}\cdot Y_{1}\ldots Y_{q}\Big)
    \label{eq:moc}
\end{align}

By Theorem \ref{thm:mocri} and Lemma \ref{lem:comb4} we obtain
\begin{align*}
  &\frac{1}{N_{1}^{\frac{3}{4}\,K}\,N_{2}^{\frac{3}{4}\,Q}} w_{K,N_{1}}^{0}(k)\; w_{Q,N_{2}}^{0}(q)\ \IE\Big(X_{1}\ldots X_{k}\cdot Y_{1}\ldots Y_{q}\Big)
  ~\leq  ~C\; N^{-\frac{K-k}{4}-\frac{Q-q}{4}}
\end{align*}

Consequently only the term with $k=K$ and $q=Q$ in \eqref{eq:moc} does not vanish in the large-$N$-limit. Thus
\begin{align*}
    & \IE\Big(\big(\frac{1}{N_{1}^{\frac{3}{4}}}\sum_{i=1}^{N_{1}}X_{i}\big)^{K}\;\big(\frac{1}{N_{2}^{\frac{3}{4}}}\sum_{j=1}^{N_{2}}Y_{j}\big)^{Q}\Big)
    ~\approx~\alpha_{1}^{\frac{K}{4}}\alpha_{2}^{\frac{Q}{4}} N^{\frac{K}{4}+\frac{Q}{4}}\IE ( X_1 \cdot ... \cdot X_K \cdot Y_1 \cdot ... \cdot Y_Q)\notag\\
    \approx~& \alpha_{1}^{\frac{K}{4}}\alpha_{2}^{\frac{Q}{4}} \left[\frac{12}{\alpha_1 (L_2 -\alpha_2)^2 + \alpha_2(L_1 -\alpha_1)^2} \right]^\frac{K+Q}{4}\notag
\cdot (L_1 - \alpha_1)^{Q/2} (L_2 - \alpha_2)^{K/2} \frac{\Gamma(\frac{K+Q+1}{4})}{\Gamma(\frac{1}{4})}
\end{align*}

\end{document}